\documentclass[11pt]{article}

\usepackage[utf8]{inputenc}
\usepackage[english]{babel}
\usepackage{enumerate}
\usepackage{amsthm,amsmath}
\usepackage{xcolor}
\usepackage{amssymb}
\usepackage{tikz}
\usepackage{aliascnt} 

\usepackage[bookmarks]{hyperref}
\hypersetup{
	pdftitle={Encoding rational Delta-triples on elliptic curves},
	pdfauthor={V.N. Schvollner},
	colorlinks=true,
	linkcolor=blue,
	citecolor=red,
	filecolor=cyan,
	urlcolor=magenta
}

\usepackage{cleveref}
\newtheorem{theorem}{Theorem}[section]
\newaliascnt{corollary}{theorem}
\newtheorem{corollary}[corollary]{Corollary}
\aliascntresetthe{corollary}

\newaliascnt{lemma}{theorem}
\newtheorem{lemma}[lemma]{Lemma}
\aliascntresetthe{lemma}

\newaliascnt{problem}{theorem}
\newtheorem{problem}[problem]{Problem}
\aliascntresetthe{problem}

\newaliascnt{proposition}{theorem}
\newtheorem{proposition}[proposition]{Proposition}
\aliascntresetthe{proposition}

\newtheorem{thmintro}{Theorem}

\theoremstyle{definition}
\newaliascnt{ejemplo}{theorem}
\newtheorem{ejemplo}[ejemplo]{Example}
\aliascntresetthe{ejemplo}

\newaliascnt{remark}{theorem}
\newtheorem{remark}[remark]{Remark}
\aliascntresetthe{remark}

\crefname{problem}{Problem}{Problems}\Crefname{problem}{Problem}{Problems}
\crefname{ejemplo}{Example}{Examples}\Crefname{ejemplo}{Example}{Examples}
\crefname{thmintro}{Theorem}{Theorems}\Crefname{thmintro}{Theorem}{Theorems}

\title{Encoding rational $\Delta$-triples on elliptic curves}
\author{Victor N. Schv\"ollner\thanks{IMASL, CONICET, San Luis, Argentina. E-mail: \texttt{vnsi9m6@gmail.com}}}
\date{}

\begin{document}
	\maketitle

	\begin{abstract}
		A \emph{$\Delta$-triple} for a natural number $n$ is a triple $(x,y,z)$ of divisors of $n$ with $1<x<y\leq z<\sqrt n$ and $n/x-x=(n/y-y)+(n/z-z)$. The numbers that admit such a triple are those with the $\Delta$ property, which governs the existence of certain $n$-simple triangles in factor graphs of split graphs. We allow the three components to be rational numbers, and show that, once the quotient $z/y=a/b$ is fixed, the rational $\Delta$-triples for $n$ are exactly the rational points of the elliptic curve $Y^2=n\,ab\,w(w+ab)\bigl(w+(a+b)^2\bigr)$ lying over the interval $-b(a+b)<w<-ab$. Since a rational $\Delta$-triple for $n$ amounts to $\Delta$-triples for infinitely many of the numbers $\kappa^2n$, the chord--tangent construction becomes a way of producing new $\Delta$-triples from old ones.

		For perfect squares, the curve always carries a rational point coming from a degenerate triple. This point has order $4$ if $a=b$, and infinite order if $a>b$. As a consequence, for every quotient $a/b>1$ there are infinitely many pairwise non-proportional $\Delta$-triples of perfect squares with $z/y=a/b$, in contrast with the fact that perfect squares have no $\Delta$-triple with $y=z$. For arbitrary $n$, there are infinitely many rational $\Delta$-triples for $n$ with a given quotient if and only if the corresponding curve has positive rank.
	\end{abstract}

	\textbf{Keywords:} $\Delta$ property; $\Delta$-triple; elliptic curve; chord--tangent construction; Nagell--Lutz theorem; perfect square; split graph \\
	\textbf{MSC 2020:} Primary 11G05; Secondary 11D25, 11A25, 05C75

	\section{Introduction}
	\label{sec:intro}

	For a split graph $S$, the combinatorics of $2$-switches on $S$ is encoded by its factor graph $\Phi(S)$, a multigraph introduced in \cite{vnsigma.factor.graph} whose induced cycles have length at most $4$ \cite{vnsigma.factor.graph.paths}. In \cite{vnsigma.delta.property} the following question was raised: for which $n$ is there a split graph whose factor graph contains an \emph{$n$-simple triangle of type~$0$}, that is, a triangle all of whose edges have multiplicity $n$ and whose three vertices have pairwise distinct degrees in $S$? The answer is governed by an arithmetic condition on $n$, the $\Delta$ property \cite[Theorem~1.2]{vnsigma.delta.property}. Writing $D_n$ for the set of positive divisors of $n$ and
	\[
	D^*_n=\bigl\{\,|a-b|:a,b\in D_n,\ ab=n\,\bigr\},\qquad
	D^+_n=\bigl\{\,x+y:x,y\in D^*_n-\{0\}\,\bigr\},
	\]
	the number $n$ has the \emph{$\Delta$ property}, written $n\in\mathbb{N}(\Delta)$, when $D^*_n\cap D^+_n\neq\varnothing$. Equivalently (\Cref{tool:caracterizacion}), $n$ admits a \emph{$\Delta$-triple}: divisors $1<x<y\leq z<\sqrt{n}$ of $n$ with
	\[
	\frac{n}{x}-x=\frac{n}{y}-y+\frac{n}{z}-z .
	\]

	\emph{Factor-difference sets.} The set $D^*_n$ is the \emph{factor-difference set} of $n$, introduced by Erd\H{o}s and Rosenfeld \cite{er-factor-difference}. They asked whether, for every $k\geq1$, there are integers $L_1<\dots<L_k$ such that $\bigl|\bigcap_{i=1}^{k}D^*_{L_i}\bigr|\geq k$ (Problem~885 in \cite{erdosproblems}). The answer is positive for $k=2$ \cite{er-factor-difference}, $k=3$ \cite{jimenez-urroz} and $k=4$ \cite{bremner}, and the question remains open in general. The $\Delta$ property is of a different nature, since it concerns the additive structure of a single factor-difference set. As the elements of $D^+_n$ are positive, $n\in\mathbb{N}(\Delta)$ if and only if $D^*_n-\{0\}$ is not \emph{sum-free}, that is, some element of $D^*_n-\{0\}$ is the sum of two, not necessarily distinct, elements of $D^*_n-\{0\}$; see \cite{tao-vu} for a survey on sum-free sets.

	The identity $(n/d+d)^2-(n/d-d)^2=4n$, which underlies Fermat's factorization method \cite{pomerance}, shows that an integer $h\geq0$ lies in $D^*_n$ if and only if $h^2+4n$ is a perfect square: if $h^2+4n=v^2$ with $v\in\mathbb{N}$, then $v\equiv h\pmod2$, and $(v-h)/2$ and $(v+h)/2$ are complementary divisors of $n$ whose difference is $h$ \cite{gs-squares}. In particular, $0\in D^*_n$ exactly when $n$ is a perfect square. This is why $0$ is excluded in the definition of $D^+_n$: otherwise the trivial relation $h=h+0$ would give every perfect square the $\Delta$ property. For $n=N^2$, the positive elements of $D^*_{N^2}$ are the legs $h$ of the Pythagorean triangles whose other leg is $2N$, so the three components of a $\Delta$-triple for $N^2$ give three such triangles, with legs $h_1=N^2/y-y$, $h_2=N^2/z-z$ and $h_3=N^2/x-x$, and the $\Delta$-triple equation reads $h_1+h_2=h_3$. For instance, the $\Delta$-triple $(9,10,25)$ for $30^2$ gives the triangles $(11,60,61)$, $(80,60,100)$ and $(91,60,109)$, with $11+80=91$ \cite{gs-squares}.

	The set $\mathbb{N}(\Delta)$ is closed under multiplication by squares: if $(x,y,z)$ is a $\Delta$-triple for $n$, then $(\kappa x,\kappa y,\kappa z)$ is a $\Delta$-triple for $\kappa^2n$ (\Cref{tool:rescaling}). Perfect squares play a special role in $\mathbb{N}(\Delta)$. On the one hand, they have no \emph{duplicated} $\Delta$-triple, that is, no $\Delta$-triple with $y=z$ \cite[Theorem~2.6]{gs-squares}. On the other hand, the question of whether infinitely many of them are \emph{$\Delta$-primitive}, that is, not of the form $\alpha^2m$ with $\alpha\geq2$ and $m\in\mathbb{N}(\Delta)$, raised in \cite{vnsigma.delta.property}, has been answered only recently \cite{gs-squares}.

	In this paper we allow the three components of a $\Delta$-triple to be rational numbers. This has two advantages.
	\begin{itemize}
		\item Rational triples for $n$ are the same thing as ordinary $\Delta$-triples for the numbers $\kappa^2n$, $\kappa\in\mathbb{N}$, divided by $\kappa$ (\Cref{lem:ell-rational-integer}). Since $\mathbb{N}(\Delta)$ is closed under multiplication by squares, this is a natural point of view: one rational triple for $n$ describes at once $\Delta$-triples for infinitely many of the numbers $n,4n,9n, 16n\dots$
		\item Once the quotient $z/y$ is fixed, the rational triples for $n$ are exactly the rational points of an elliptic curve lying on a certain arc of that curve. On an elliptic curve, new rational points are produced from old ones by drawing lines (the chord--tangent construction). Hence new $\Delta$-triples can be produced from old ones.
	\end{itemize}
	The second item is our first result. Here a \emph{rational $\Delta$-triple} for $n$ is a triple of rational numbers $0<x<y\leq z<\sqrt n$ satisfying the equation above (\Cref{sec:rational}).

	\begin{thmintro}[\Cref{prop:ell-curve}]
		\label{thm:A}
		Let $n\in\mathbb{N}$, let $a\geq b\geq1$ be coprime integers, and put $s=a+b$. The map
		\[
		(x,y,z)\longmapsto\Bigl(-\frac{bs\,x}{y},\ \frac{ab s\,x\,(sy-bx)}{y}\Bigr)
		\]
		is a bijection from the set of rational $\Delta$-triples $(x,y,z)$ for $n$ with $z/y=a/b$ onto the set of rational points $(w,Y)$ of the elliptic curve
		\[
		E^{(n)}_{a,b}:\qquad Y^2=n\,ab\;w\,(w+ab)\,\bigl(w+s^2\bigr)
		\]
		with $-bs<w<-ab$ and $Y>0$.
	\end{thmintro}

	The interval $-bs<w<-ab$ is the \emph{window}. The curve depends on $n$ only up to square factors, so a single curve $E_{a,b}=E^{(1)}_{a,b}$ describes the $\Delta$-triples with $z/y=a/b$ of all perfect squares at once (\Cref{rem:ell-scaling}). Moreover, $E_{a,b}$ always contains the rational point
	\[
	P=\bigl(-bs,\ ab^2s\bigr),
	\]
	which sits at the left end of the window and corresponds to the degenerate triple $(b/a,\,b/a,\,1)$ (\Cref{rem:ell-ends}). In terms of factor-difference sets, this triple is the trivial relation $h=h+0$: indeed $x=y$, and $1/z-z=0$. Thus $P$ comes from the element $0$, which only the factor-difference sets of perfect squares contain. Everything then depends on the order of $P$. The duplication formula and the theorem of Nagell and Lutz show that $P$ has order $4$ if $a=b$, and infinite order if $a>b$ (\Cref{lem:ell-order-P}). In the second case, an argument based on the structure of the real points of the curve and on Kronecker's approximation theorem shows that infinitely many odd multiples of $P$ lie over the window (\Cref{lem:ell-density}). This is our main result.

	\begin{thmintro}[\Cref{thm:ell-every-quotient}]
		\label{thm:B}
		Let $a>b\geq1$ be coprime integers. Then there are infinitely many pairwise non-proportional $\Delta$-triples $(x,y,z)$ of perfect squares with $z/y=a/b$.
	\end{thmintro}

	Thus every quotient $a/b>1$ occurs among the $\Delta$-triples of perfect squares, while the quotient $1$ never occurs \cite[Theorem~2.6]{gs-squares}; on the curves, both facts are explained by the order of the single point $P$. In particular, for $a>b$ the chord--tangent construction turns the trivial relation $h=h+0$ into infinitely many genuine relations $h_1+h_2=h_3$ among positive elements of factor-difference sets of perfect squares. For $a=b$, the absence of rational points of $E_{1,1}$ over the window gives an alternative proof of the fact that perfect squares have no duplicated $\Delta$-triples, based on the rank of $E_{1,1}$ instead of the theorem of Euler on products of four terms in arithmetic progression used in \cite{gs-squares} (\Cref{rem:ell-duplicated}). The triples produced by \Cref{thm:B} can be very large: the least perfect square with a $\Delta$-triple with $z/y=3$ is $18837^2$ (\Cref{sec:conclusion}), while the first odd multiple of $P$ over the window of $E_{3,1}$ gives one for $N^2$ with $N=488650465808508899094963$ (\Cref{ex:ell-squares}).

	For $n$ not a perfect square, the degenerate point is not rational, and there is no free point to start from. The question becomes whether the curve has rational points of infinite order.

	\begin{thmintro}[\Cref{prop:ell-rank}]
		\label{thm:C}
		Let $n\in\mathbb{N}$ and let $a\geq b\geq1$ be coprime integers. There are infinitely many rational $\Delta$-triples for $n$ with $z/y=a/b$ if and only if $E^{(n)}_{a,b}$ has positive rank. In that case, there are infinitely many pairwise non-proportional $\Delta$-triples with $z/y=a/b$ among the numbers $\kappa^2n$, $\kappa\in\mathbb{N}$.
	\end{thmintro}

	In particular, a single $\Delta$-triple whose point on the curve has infinite order produces infinitely many non-proportional ones (\Cref{cor:ell-triple-infinite-order}). With PARI/GP \cite{pari} we checked that this is the case for every $\Delta$-triple of every $n\leq10^5$: there are $4058$ such triples, for $1780$ values of $n$.

	The curves $E^{(n)}_{a,b}$ see a number together with all its multiples by squares, so they cannot detect $\Delta$-primitivity (\Cref{subsec:ell-limits}). In this sense, the present approach is complementary to the arithmetic one of \cite{gs-squares}.

	The paper is meant to be readable without any background on elliptic curves. Besides elementary algebra, we use only four classical facts on elliptic curves, quoted with references where they are needed: the group law (see {\cite[Chapter~I]{silverman-tate}}), the theorem of Nagell and Lutz (\Cref{fact:ell-NL}), the structure of the real points (\Cref{fact:ell-real}), and Mordell's theorem (only in \Cref{prop:ell-rank}), together with Kronecker's approximation theorem (\Cref{fact:kronecker}).

	\emph{Organization.} \Cref{sec:prelim} collects the facts on the $\Delta$ property that we need. \Cref{sec:rational} introduces rational $\Delta$-triples and shows that they are rescaled $\Delta$-triples (\Cref{lem:ell-rational-integer}). In \Cref{subsec:ell-curve} we attach to each quotient $a/b$ the curve $E^{(n)}_{a,b}$ and prove \Cref{thm:A}. \Cref{subsec:ell-chords} explains the chord--tangent construction, and \Cref{subsec:ell-real} studies the real points of the curve. Perfect squares are treated in \Cref{subsec:ell-squares}, where \Cref{thm:B} is proved, and arbitrary $n$ in \Cref{subsec:ell-rank}, where \Cref{thm:C} is proved. \Cref{subsec:ell-limits} discusses the limits of the method, and \Cref{sec:conclusion} contains concluding remarks and open problems.

	\section{Preliminaries}
	\label{sec:prelim}

	We use the notation of \cite{vnsigma.delta.property}.

	\begin{proposition}[{\cite[Proposition~2.3]{vnsigma.delta.property}}]
		\label{tool:caracterizacion}
		A natural number $n$ lies in $\mathbb{N}(\Delta)$ if and only if there exists $(x,y,z)\in D_n^3$ with $1<x<y\leq z<\sqrt{n}$ and
		\begin{equation}
			\frac{n}{x}-x \;=\; \frac{n}{y}-y + \frac{n}{z}-z.
			\label{eq:tool-delta}
		\end{equation}
		Such a triple is a \emph{$\Delta$-triple} for $n$. Equation \eqref{eq:tool-delta} is equivalent to
		\begin{equation}
			(xy+xz-yz)\,n \;=\; xyz\,(z+y-x).
			\label{eq:tool-delta2}
		\end{equation}
	\end{proposition}

	A $\Delta$-triple $(x,y,z)$ is \emph{duplicated} if $y=z$, and \emph{generic} if $y<z$ \cite{vnsigma.delta.property,gs-squares}. Perfect squares have no duplicated $\Delta$-triple \cite[Theorem~2.6]{gs-squares}.

	\begin{proposition}[{\cite[Proposition~3.1]{vnsigma.delta.property}}]
		\label{tool:rescaling}
		If $(x,y,z)$ is a $\Delta$-triple for $n$, then $(\alpha x,\alpha y,\alpha z)$ is a $\Delta$-triple for $\alpha^2n$ for every $\alpha\in\mathbb{N}$. In particular, if $n\in\mathbb{N}(\Delta)$, then $\alpha^2 n\in\mathbb{N}(\Delta)$ for all $\alpha\in\mathbb{N}$.
	\end{proposition}

	An element $n\in\mathbb{N}(\Delta)$ is \emph{$\Delta$-primitive} if there are no integers $\alpha,m\geq2$ with $n=\alpha^2m$ and $m\in\mathbb{N}(\Delta)$. There are infinitely many $\Delta$-primitive numbers \cite[Theorem~3.10]{vnsigma.delta.property}, and it was conjectured in \cite[Conjecture~3.11]{vnsigma.delta.property} that infinitely many of them are perfect squares. This conjecture was proved in \cite{gs-squares}, using a polynomial family of $\Delta$-triples whose square values are governed by the Pell equation $u^2-3v^2=1$. Following \cite{gs-squares}, put
	\[
	\mathcal{S}=\{\,N\in\mathbb{N}:\ N^2\in\mathbb{N}(\Delta)\,\}.
	\]
	By \Cref{tool:rescaling}, if $N\in\mathcal{S}$ and $N\mid M$, then $M\in\mathcal{S}$. Moreover, $N^2$ is $\Delta$-primitive if and only if $N\in\mathcal{S}$ and no divisor $d<N$ of $N$ lies in $\mathcal{S}$ \cite[Proposition~3.3]{gs-squares}. For instance, $30^2$ is the least $\Delta$-primitive square, and $(9,10,25)$ is its only $\Delta$-triple \cite[Example~3.5]{gs-squares}.

\section{Rational $\Delta$-triples}
\label{sec:rational}

Let $n\in\mathbb{N}$. A \emph{rational $\Delta$-triple} for $n$ is a triple $(x,y,z)$ of rational numbers with $0<x<y\leq z<\sqrt{n}$ that satisfies \eqref{eq:tool-delta}. We write $\mathcal{R}_n$ for the set of rational $\Delta$-triples for $n$ and, for coprime integers $a\geq b\geq1$,
\[
\mathcal{R}_n(a/b)=\bigl\{(x,y,z)\in\mathcal{R}_n:\ z/y=a/b\bigr\}.
\]
Note that neither divisibility nor the condition $x>1$ is required. Every $\Delta$-triple for $n$ is a rational $\Delta$-triple for $n$. Since $z\geq y$, the quotient $z/y$ is a rational number $\geq1$, and we always write it in lowest terms as $a/b$ with $a\geq b\geq1$. So every element of $\mathcal{R}_n$ lies in exactly one set $\mathcal{R}_n(a/b)$, and $\mathcal{R}_n(1)$ consists of the triples with $y=z$.

Multiplying \eqref{eq:tool-delta} by $xyz$ shows that \eqref{eq:tool-delta2} also holds for rational $\Delta$-triples. Its right-hand side $xyz(y+z-x)$ is positive, because $x,y,z>0$ and $y>x$. Hence $xy+xz-yz>0$, and
\begin{equation}
	\label{eq:ell-n-from-triple}
	n=\frac{xyz\,(y+z-x)}{xy+xz-yz}.
\end{equation}
Thus a rational $\Delta$-triple determines $n$. Moreover, replacing $(x,y,z)$ by $(cx,cy,cz)$ with $c>0$ multiplies the right-hand side of \eqref{eq:ell-n-from-triple} by $c^2$.

\begin{ejemplo}
	\label{ex:ell-24}
	The triple $(1,\tfrac43,\tfrac83)$ lies in $\mathcal{R}_{24}(2)$. Indeed $0<1<\tfrac43<\tfrac83<\sqrt{24}$, $z/y=2$, and
	\[
	\frac{24}{1}-1=23,\qquad \frac{24}{4/3}-\frac43=\frac{50}{3},\qquad \frac{24}{8/3}-\frac83=\frac{19}{3},\qquad \frac{50}{3}+\frac{19}{3}=23 .
	\]
\end{ejemplo}

The next lemma says that rational $\Delta$-triples are just rescaled $\Delta$-triples.

\begin{lemma}
	\label{lem:ell-rational-integer}
	Let $n,\kappa\in\mathbb{N}$.
	\begin{enumerate}
		\item The map $(x,y,z)\mapsto(\kappa x,\kappa y,\kappa z)$ is a bijection from $\mathcal{R}_n$ onto $\mathcal{R}_{\kappa^2n}$, and it preserves the quotient $z/y$.
		\item If $(x',y',z')$ is a $\Delta$-triple for $\kappa^2n$, then $(x'/\kappa,y'/\kappa,z'/\kappa)\in\mathcal{R}_n$.
		\item Let $(x,y,z)\in\mathcal{R}_n$, and let $L$ be the least common multiple of the denominators of the six rational numbers
		\[
		x,\quad y,\quad z,\quad n/x,\quad n/y,\quad n/z .
		\]
		If $\kappa$ is a multiple of $L$ with $\kappa x>1$, then $(\kappa x,\kappa y,\kappa z)$ is a $\Delta$-triple for $\kappa^2n$.
	\end{enumerate}
\end{lemma}

\begin{proof}
	For every rational $t>0$ we have
	\begin{equation}
		\label{eq:ell-scaling}
		\frac{\kappa^2n}{\kappa t}-\kappa t=\kappa\Bigl(\frac nt-t\Bigr).
	\end{equation}
	
	(1) The inequalities $0<x<y\leq z<\sqrt n$ are equivalent to $0<\kappa x<\kappa y\leq\kappa z<\kappa\sqrt n=\sqrt{\kappa^2n}$. By \eqref{eq:ell-scaling}, equation \eqref{eq:tool-delta} for $(\kappa x,\kappa y,\kappa z)$ and $\kappa^2n$ is equation \eqref{eq:tool-delta} for $(x,y,z)$ and $n$, multiplied by $\kappa$. Hence $(x,y,z)\in\mathcal{R}_n$ if and only if $(\kappa x,\kappa y,\kappa z)\in\mathcal{R}_{\kappa^2n}$. The inverse map is division by $\kappa$. Clearly $(\kappa z)/(\kappa y)=z/y$.
	
	(2) A $\Delta$-triple for $\kappa^2n$ lies in $\mathcal{R}_{\kappa^2n}$. Apply the inverse map of (1).
	
	(3) By (1), $(\kappa x,\kappa y,\kappa z)\in\mathcal{R}_{\kappa^2n}$. Since $L\mid\kappa$, the numbers $\kappa x$, $\kappa y$, $\kappa z$ are integers, and so is $\kappa^2n/(\kappa x)=\kappa\cdot(n/x)$; that is, $\kappa x$ divides $\kappa^2n$. The same holds for $\kappa y$ and $\kappa z$. Together with $\kappa x>1$, \Cref{tool:caracterizacion} shows that $(\kappa x,\kappa y,\kappa z)$ is a $\Delta$-triple for $\kappa^2n$.
\end{proof}

\begin{ejemplo}
	For $(1,\tfrac43,\tfrac83)\in\mathcal{R}_{24}$ the six numbers of \Cref{lem:ell-rational-integer}(3) are $1,\tfrac43,\tfrac83,24,18,9$, so $L=3$. With $\kappa=3$ we obtain the $\Delta$-triple $(3,4,8)$ for $216=3^2\cdot24$: indeed $72-3=69=(54-4)+(27-8)$.
\end{ejemplo}

Two triples $T$ and $T'$ of positive rational numbers are \emph{proportional} if $T'=c\,T$ for some rational $c>0$. Proportional triples have the same quotient $z/y$.

\begin{corollary}
	\label{cor:ell-classes}
	Let $n\in\mathbb{N}$ and let $a\geq b\geq1$ be coprime. If $\mathcal{R}_n(a/b)$ is infinite, then there are infinitely many pairwise non-proportional triples that are $\Delta$-triples for numbers of the form $\kappa^2n$, $\kappa\in\mathbb{N}$, and have $z/y=a/b$. In particular, if $\mathcal{R}_1(a/b)$ is infinite, then there are infinitely many pairwise non-proportional $\Delta$-triples of perfect squares with $z/y=a/b$.
\end{corollary}

\begin{proof}
	Choose pairwise distinct elements $t_1,t_2,\dots$ of $\mathcal{R}_n(a/b)$. By \Cref{lem:ell-rational-integer}(3), for each $i$ there exists $\kappa_i\in\mathbb{N}$ such that $T_i=\kappa_i t_i$ is a $\Delta$-triple for $\kappa_i^2n$, and $T_i$ has the same quotient $z/y=a/b$ as $t_i$.
	
	Suppose that $T_j=c\,T_i$ for some rational $c>0$. Then $t_j=\lambda\,t_i$ with $\lambda=c\kappa_i/\kappa_j>0$. Since $t_i$ and $t_j$ both lie in $\mathcal{R}_n$, applying \eqref{eq:ell-n-from-triple} to $t_i$ gives $n$, and applying it to $t_j=\lambda t_i$ gives $\lambda^2n$, because multiplying a triple by $\lambda$ multiplies the right-hand side of \eqref{eq:ell-n-from-triple} by $\lambda^2$. Hence $\lambda=1$, so $t_j=t_i$ and $j=i$. Thus the $T_i$ are pairwise non-proportional.
	
	For the last assertion take $n=1$: the numbers $\kappa^2n$ are then exactly the perfect squares.
\end{proof}

\section{The curve attached to a quotient $z/y=a/b$}
\label{subsec:ell-curve}

Throughout the rest of the paper we fix $n\in\mathbb{N}$ and coprime integers $a\geq b\geq1$, and we write
\[
s=a+b .
\]
Consider the curve
\begin{equation}
	\label{eq:ell-curve}
	E^{(n)}_{a,b}:\qquad Y^2=n\,ab\;w\,(w+ab)\,\bigl(w+s^2\bigr),
\end{equation}
and write $E_{a,b}=E^{(1)}_{a,b}$. The cubic polynomial on the right-hand side has the three distinct roots $-s^2<-ab<0$, so \eqref{eq:ell-curve} defines an elliptic curve. A \emph{rational point} of $E^{(n)}_{a,b}$ is a pair $(w,Y)\in\mathbb{Q}^2$ satisfying \eqref{eq:ell-curve}. The interval
\[
-bs<w<-ab
\]
is called the \emph{window}. It is non-empty, because $bs=ab+b^2>ab$. We say that a point $(w,Y)$ lies \emph{over the window} if $w$ lies in the window.

\begin{theorem}
	\label{prop:ell-curve}
	Let $n\in\mathbb{N}$ and let $a\geq b\geq1$ be coprime, $s=a+b$.
	\begin{enumerate}
		\item If $(x,y,z)\in\mathcal{R}_n(a/b)$, then
		\begin{equation}
			\label{eq:ell-forward}
			w=-\frac{bs\,x}{y},\qquad Y=\frac{ab s\,x\,(sy-bx)}{y}
		\end{equation}
		is a rational point of $E^{(n)}_{a,b}$, with $w$ in the window and $Y>0$.
		\item Conversely, if $(w,Y)$ is a rational point of $E^{(n)}_{a,b}$ with $w$ in the window and $Y>0$, then
		\begin{equation}
			\label{eq:ell-backward}
			x=\frac{Y}{ab\,(w+s^2)},\qquad y=\frac{-sY}{a\,w\,(w+s^2)},\qquad z=\frac{-sY}{b\,w\,(w+s^2)}
		\end{equation}
		defines an element $(x,y,z)$ of $\mathcal{R}_n(a/b)$.
		\item The maps in (1) and (2) are inverse to each other. Hence $\mathcal{R}_n(a/b)$ is in bijection with the set of rational points $(w,Y)$ of $E^{(n)}_{a,b}$ with $w$ in the window and $Y>0$.
	\end{enumerate}
\end{theorem}

\begin{proof}
	We proceed in five steps.
	
		\emph{Step 1: a change of parameters.} Let $x,y,z$ be positive rational numbers with $z/y=a/b$, that is, $bz=ay$. Put
	\[
	\lambda=\frac yb=\frac za,\qquad r=\frac x\lambda=\frac{bx}{y} .
	\]
	Then $\lambda$ and $r$ are positive rationals, and $(x,y,z)=\lambda\,(r,b,a)$. Conversely, for any positive rationals $r$ and $\lambda$, the triple $\lambda\,(r,b,a)$ consists of positive rationals and has quotient $z/y=a/b$; and the two constructions are mutually inverse. Hence $(x,y,z)\mapsto(r,\lambda)$ is a bijection from the set of triples of positive rationals with $z/y=a/b$ onto the set of pairs of positive rationals. Note that $r=bx/y$ depends only on the quotient $x/y$, whereas $\lambda$ carries the scaling: replacing $(x,y,z)$ by $c\,(x,y,z)$, with $c>0$ rational, replaces $(r,\lambda)$ by $(r,c\lambda)$.
	
	\emph{Step 2: the equation.} For $(x,y,z)=(r\lambda,b\lambda,a\lambda)$ the right-hand side of \eqref{eq:tool-delta} is
	\[
	\frac{n}{b\lambda}+\frac{n}{a\lambda}-(a+b)\lambda=\frac{ns}{ab\,\lambda}-s\lambda,
	\]
	and the left-hand side is $\frac{n}{r\lambda}-r\lambda$. Multiplying by $\lambda$ and collecting the terms with $\lambda^2$, equation \eqref{eq:tool-delta} becomes
	\begin{equation}
		\label{eq:ell-E1}
		(s-r)\,\lambda^2=\frac{ns}{ab}-\frac nr=\frac{n\,(sr-ab)}{ab\,r}.
	\end{equation}
	
	\emph{Step 3: the inequalities.} The condition $x<y$ means $r<b$, and $y\leq z$ holds because $b\leq a$. Assume $r<b$. Then $s-r>a>0$, so the left-hand side of \eqref{eq:ell-E1} is positive, which forces $sr>ab$, that is $r>ab/s$. Then \eqref{eq:ell-E1} reads
	\begin{equation}
		\label{eq:ell-E2}
		\lambda^2=\frac{n\,(sr-ab)}{ab\,r\,(s-r)} .
	\end{equation}
	It remains to see what $z<\sqrt n$ means. By \eqref{eq:ell-E2}, and multiplying by the positive number $ab\,r(s-r)/n$, the condition $a^2\lambda^2<n$ becomes $a^2(sr-ab)<ab\,r(s-r)$. Expanding, using $s=a+b$ and dividing by $a$, this is
	\[
	br^2+(a^2-b^2)\,r-a^2b<0,\qquad\text{that is,}\qquad (br+a^2)(r-b)<0 .
	\]
	Since $br+a^2>0$, this is again the condition $r<b$. Summing up:
	\begin{equation}
		\label{eq:ell-summary}
		(r\lambda,b\lambda,a\lambda)\in\mathcal{R}_n(a/b)\iff \frac{ab}{s}<r<b\quad\text{and}\quad\lambda>0\ \text{satisfies \eqref{eq:ell-E2}}.
	\end{equation}
	(In particular the condition $z<\sqrt n$ is automatic.)
	
	\emph{Step 4: when is $\lambda$ rational?} Fix a rational $r$ with $ab/s<r<b$. The right-hand side of \eqref{eq:ell-E2} is a positive rational number. Multiplying \eqref{eq:ell-E2} by $\bigl(ab\,r(s-r)\bigr)^2$, we see that the positive rationals $\lambda$ satisfying \eqref{eq:ell-E2} correspond bijectively, via $Y'=\lambda\,ab\,r(s-r)$, to the positive rationals $Y'$ satisfying
	\[
	(Y')^2=abn(sr-ab)(s-r)r.
	\]
	
	\emph{Step 5: the change of variable $w=-sr$.} Put $w=-sr$, so that $r=-w/s$. Then
	\[
	sr-ab=-(w+ab),\qquad s-r=\frac{w+s^2}{s},
	\]
	and hence
	\[
	r\,(sr-ab)\,(s-r)=\Bigl(-\frac ws\Bigr)\cdot\bigl(-(w+ab)\bigr)\cdot\frac{w+s^2}{s}=\frac{w\,(w+ab)\,(w+s^2)}{s^2}.
	\]
	Putting $Y=sY'$, the equation of Step 4 becomes \eqref{eq:ell-curve}. The condition $ab/s<r<b$ becomes $-bs<w<-ab$, and $Y'>0$ becomes $Y>0$.
	
	\emph{The formulas.} Given $(x,y,z)$, Step 1 gives $r=bx/y$ and $\lambda=y/b$, so $w=-sr=-bsx/y$ and
	\[
	Y=s\,\lambda\,ab\,r(s-r)=s\,ab\,x\Bigl(s-\frac{bx}{y}\Bigr)=\frac{abs\,x\,(sy-bx)}{y},
	\]
	which is \eqref{eq:ell-forward}. Conversely, given $(w,Y)$, we have $r=-w/s$ and $s\,ab\,r(s-r)=-ab\,w(w+s^2)/s$, so
	\[
	\lambda=\frac{Y}{s\,ab\,r(s-r)}=\frac{-sY}{ab\,w\,(w+s^2)},
	\]
	and $x=r\lambda$, $y=b\lambda$, $z=a\lambda$ are the formulas \eqref{eq:ell-backward}. Each step above is reversible, so (1), (2) and (3) follow from \eqref{eq:ell-summary}.
\end{proof}

In short: \emph{the rational $\Delta$-triples for $n$ with $z/y=a/b$ are the rational points of $E^{(n)}_{a,b}$ over the window} (taking $Y>0$). Note that only $w$ matters: by \eqref{eq:ell-forward}, $w$ determines $x/y$, and then $Y$ (taken positive) determines the size. The two points $(w,Y)$ and $(w,-Y)$ will be regarded as giving the same triple, namely the one obtained from $(w,|Y|)$.

\begin{remark}
	\label{rem:ell-scaling}
	Strictly speaking, each number has its own curve: the number $m$ has the curve $E^{(m)}_{a,b}$. But the curve depends on $m$ only up to square factors, in the following sense. Let $n,\kappa\in\mathbb{N}$. Comparing the equations of $E^{(n)}_{a,b}$ and $E^{(\kappa^2n)}_{a,b}$, we see that $(w,Y)$ satisfies the first one if and only if $(w,\kappa Y)$ satisfies the second one. Hence
	\[
	\iota_\kappa\colon E^{(n)}_{a,b}\longrightarrow E^{(\kappa^2n)}_{a,b},\qquad \iota_\kappa(w,Y)=(w,\kappa Y),
	\]
	is a bijection; it and its inverse take rational points to rational points, and it leaves $w$ unchanged, so it preserves the window and the sign of $Y$. By \Cref{prop:ell-curve} it therefore matches $\mathcal{R}_n(a/b)$ with $\mathcal{R}_{\kappa^2n}(a/b)$, in accordance with \Cref{lem:ell-rational-integer}(1): by \eqref{eq:ell-forward}, multiplying $(x,y,z)$ by $\kappa$ leaves $w$ unchanged and multiplies $Y$ by $\kappa$.
	
	In particular, for every $\kappa\in\mathbb{N}$ the curve of the perfect square $\kappa^2$, namely $E^{(\kappa^2)}_{a,b}$, is identified by $\iota_\kappa$ with the single curve $E_{a,b}=E^{(1)}_{a,b}$. This is the first reason why perfect squares are special in \Cref{subsec:ell-squares}: one curve describes the $\Delta$-triples of all of them at once. The second reason is that, by \Cref{rem:ell-ends}, the point at the left end of the window is rational exactly when the number is a perfect square; on $E_{a,b}$ it is the point $P$ used in \Cref{subsec:ell-squares}.
\end{remark}

\begin{remark}[The two ends of the window]
	\label{rem:ell-ends}
	At the left end $w=-bs$ we have $r=b$, that is $x=y$, and \eqref{eq:ell-E2} gives $\lambda^2=n/a^2$, so $z=a\lambda=\sqrt n$. So the left end corresponds to the \emph{degenerate triple}
	\[
	\Bigl(\frac{b\sqrt n}{a},\ \frac{b\sqrt n}{a},\ \sqrt n\Bigr),
	\]
	which satisfies \eqref{eq:tool-delta} for a trivial reason: $x=y$, and $n/z-z=0$. It violates both strict inequalities $x<y$ and $z<\sqrt n$, so it is not a rational $\Delta$-triple. The corresponding point of the curve is $(-bs,\,ab^2s\sqrt n)$, since
	\[
	n\,ab\cdot(-bs)\,(-bs+ab)\,(-bs+s^2)=n\,ab\cdot(-bs)(-b^2)(as)=n\,a^2b^4s^2 .
	\]
	It is a rational point exactly when $n$ is a perfect square. This point is the key to \Cref{subsec:ell-squares}. At the right end $w=-ab$ we have $Y=0$ and $\lambda=0$: the triple collapses to $(0,0,0)$.
\end{remark}

\begin{ejemplo}
	\label{ex:ell-24-curve}
	Let $n=24$ and $a/b=2$, so $s=3$. The curve and the window are
	\[
	E^{(24)}_{2,1}:\ Y^2=48\,w\,(w+2)\,(w+9),\qquad -3<w<-2 .
	\]
	For the triple $(1,\tfrac43,\tfrac83)$ of \Cref{ex:ell-24}, \eqref{eq:ell-forward} gives
	\[
	w=-\frac{3\cdot1}{4/3}=-\frac94,\qquad Y=\frac{2\cdot3\cdot1\cdot(4-1)}{4/3}=\frac{27}{2},
	\]
	and indeed $48\cdot(-\tfrac94)(-\tfrac14)(\tfrac{27}{4})=\tfrac{729}{4}=\bigl(\tfrac{27}2\bigr)^2$, with $-3<-\tfrac94<-2$. We call this point $P_0=(-\tfrac94,\tfrac{27}2)$.
\end{ejemplo}

\section{Producing new points: chords and tangents}
\label{subsec:ell-chords}

From now on, $E=E^{(n)}_{a,b}$ and $f(w)=n\,ab\,w(w+ab)(w+s^2)$, so that $E$ is $Y^2=f(w)$. The basic observation is that a line meets $E$ in three points, and if two of them are rational, so is the third.

\begin{lemma}[Chord--tangent construction; see {\cite[Chapter~I]{silverman-tate}}]
	\label{lem:ell-third-point}
	Let $P_1=(w_1,Y_1)$ and $P_2=(w_2,Y_2)$ be points of $E$ with real coordinates, and assume that either $w_1\neq w_2$, or $P_1=P_2$ and $Y_1\neq0$. Let $Y=mw+\beta$ be the line through $P_1$ and $P_2$ (the tangent line to $E$ at $P_1$ if $P_1=P_2$), so that
	\[
	m=\frac{Y_2-Y_1}{w_2-w_1}\quad(w_1\neq w_2),\qquad m=\frac{f'(w_1)}{2Y_1}\quad(P_1=P_2).
	\]
	Then the line meets $E$ in a third point $P_3=(w_3,\,mw_3+\beta)$, where
	\begin{equation}
		\label{eq:ell-vieta}
		w_3=\frac{m^2}{n\,ab}-ab-s^2-w_1-w_2 .
	\end{equation}
	More precisely, $w_1,w_2,w_3$ are the three roots, counted with multiplicity, of $g(w)=f(w)-(mw+\beta)^2$. If $P_1$ and $P_2$ are rational points, then $m$, $\beta$ and $P_3$ are rational.
\end{lemma}

This construction turns the set of rational points into a group. We add to $E$ one \emph{point at infinity} $O$, through which all vertical lines pass, and we define, for points $P=(w,Y)$, $P_1$, $P_2$:
\begin{itemize}
	\item $-P=(w,-Y)$, the mirror image of $P$, and $-O=O$;
	\item $P_1+P_2=-P_3$, where $P_3$ is the third point of \Cref{lem:ell-third-point};
	\item $P+(-P)=O$ (the vertical line through $P$ and $-P$ ``meets $E$ again at $O$''), and $P+O=O+P=P$.
\end{itemize}
Since $-P$ and $P$ have the same $w$, \eqref{eq:ell-vieta} gives $w(P_1+P_2)=w_3$. With this addition, the rational points of $E$ together with $O$ form an abelian group $E(\mathbb{Q})$ with neutral element $O$, in which the opposite of $P$ is $-P$. Likewise, the real points of $E$ together with $O$ form an abelian group $E(\mathbb{R})$, which contains $E(\mathbb{Q})$ as a subgroup (see \cite[Chapter~I]{silverman-tate}).

In \cite{silverman-tate} curves are written as $V^2=u^3+Au^2+Bu+C$. Our curve has this form after the change of variables
\begin{equation}
	\label{eq:ell-weierstrass}
	u=n\,ab\,w,\qquad V=n\,ab\,Y,\qquad\qquad V^2=u\,\bigl(u+n\,a^2b^2\bigr)\bigl(u+n\,ab\,s^2\bigr),
\end{equation}
which follows by multiplying \eqref{eq:ell-curve} by $(n\,ab)^2$. This change of variables maps lines to lines, so it does not affect the construction above. What it achieves is that the cubic becomes \emph{monic}: the equation of $E^{(n)}_{a,b}$ already has integer coefficients, but its leading coefficient is $n\,ab$, whereas \cite{silverman-tate} and \Cref{fact:ell-NL} below require it to be $1$. Expanding, \eqref{eq:ell-weierstrass} reads $V^2=u^3+Au^2+Bu$ with the integers $A=n\,ab\,(ab+s^2)$ and $B=n^2a^3b^3s^2$. We call \eqref{eq:ell-weierstrass} the \emph{Weierstrass model} of $E$.

For $P\in E(\mathbb{Q})$ and $j\in\mathbb{N}$ we write $jP=P+\cdots+P$ ($j$ times). The point $P$ has \emph{finite order} if $jP=O$ for some $j\geq1$; otherwise $P$ has \emph{infinite order}, and then the points $P,2P,3P,\dots$ are pairwise distinct. If $P$ has finite order, so does every multiple $jP$.

\begin{theorem}[Nagell--Lutz; see {\cite[Chapter~II]{silverman-tate}}]
	\label{fact:ell-NL}
	Let $A,B,C\in\mathbb{Z}$. If $(u,V)$ is a rational point of finite order on $V^2=u^3+Au^2+Bu+C$, then $u$ and $V$ are integers.
\end{theorem}

We shall use \Cref{fact:ell-NL} in the following way: \emph{if some multiple $jP$ has a non-integral $u$-coordinate in the Weierstrass model \eqref{eq:ell-weierstrass}, then $P$ has infinite order.} To compute doublings we use the following formula.

\begin{lemma}[Duplication formula; see {\cite[Chapter~I]{silverman-tate}}]
	\label{lem:ell-duplication}
	Let $A,B\in\mathbb{Q}$ and let $Q=(u_0,V_0)$ be a point of $V^2=u^3+Au^2+Bu$ with $V_0\neq0$. Then
	\[
	u(2Q)=\frac{(u_0^2-B)^2}{4V_0^2}.
	\]
\end{lemma}

\begin{ejemplo}
	\label{ex:ell-24-chords}
	We continue \Cref{ex:ell-24-curve}: $f(w)=48(w^3+11w^2+18w)$ and $P_0=(-\tfrac94,\tfrac{27}2)$.
	\begin{itemize}
		\item \emph{Tangent at $P_0$.} We have $f'(w)=48(3w^2+22w+18)$ and $f'(-\tfrac94)=48\cdot(-\tfrac{261}{16})=-783$, so $m=-783/27=-29$. By \eqref{eq:ell-vieta},
		\[
		w(2P_0)=\frac{841}{48}-2-9+\frac94+\frac94=\frac{529}{48},\qquad 2P_0=\Bigl(\frac{529}{48},\ \frac{17825}{48}\Bigr).
		\]
		Here $w>0$ is outside the window, so $2P_0$ gives no triple.
		\item \emph{Chord through $P_0$ and $2P_0$.} The slope is $m=\tfrac{17177}{637}$, and \eqref{eq:ell-vieta} gives
		\[
		w(3P_0)=-\frac{7502121}{1623076}\approx-4.622 .
		\]
		This is not in the window either.
		\item \emph{Order of $P_0$.} The Weierstrass model \eqref{eq:ell-weierstrass} is $V^2=u(u+96)(u+432)$ with $u=48w$. Since $u(3P_0)=48\,w(3P_0)=-\tfrac{90025452}{405769}$ is not an integer, $3P_0$ has infinite order by \Cref{fact:ell-NL}, and hence so does $P_0$.
	\end{itemize}
\end{ejemplo}

\section{The oval and the unbounded branch}
\label{subsec:ell-real}

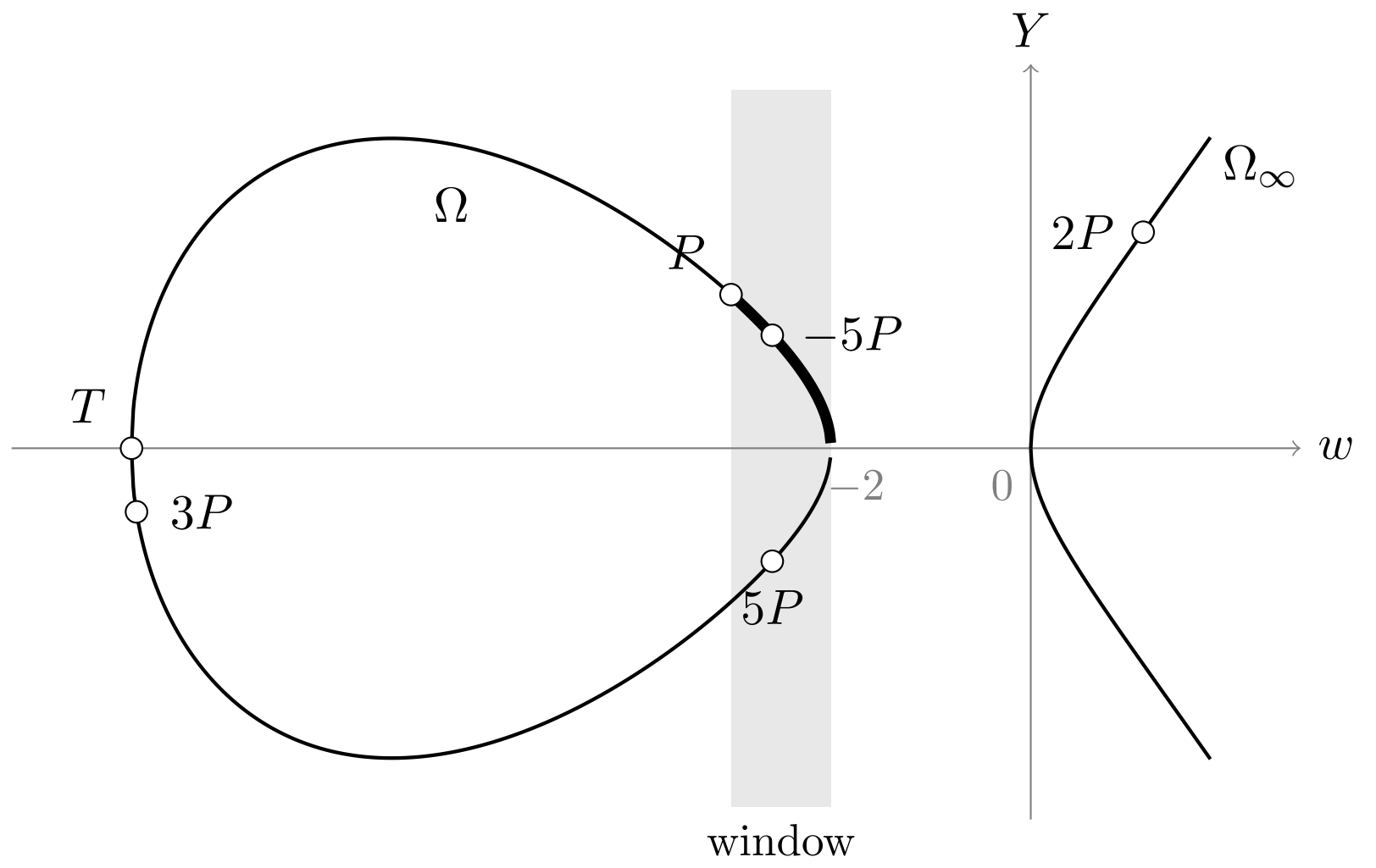
\begin{figure}[ht]
	\centering
	\begin{tikzpicture}[xscale=0.78,yscale=0.2,
		pt/.style={circle,draw=black,fill=white,inner sep=1.7pt}]
		\fill[gray!18] (-3,-14) rectangle (-2,14);
		\draw[->,gray] (-10.2,0) -- (2.7,0) node[right,black] {$w$};
		\draw[->,gray] (0,-14.5) -- (0,15) node[above,black] {$Y$};
		\draw[thick,domain=-9:-2,samples=400,smooth] plot (\x,{sqrt(max(0,2*\x*(\x+2)*(\x+9)))});
		\draw[thick,domain=-9:-2,samples=400,smooth] plot (\x,{-sqrt(max(0,2*\x*(\x+2)*(\x+9)))});
		\draw[thick,domain=0:1.8,samples=200,smooth] plot (\x,{sqrt(max(0,2*\x*(\x+2)*(\x+9)))});
		\draw[thick,domain=0:1.8,samples=200,smooth] plot (\x,{-sqrt(max(0,2*\x*(\x+2)*(\x+9)))});
		\draw[line width=2.4pt,domain=-3:-2,samples=100,smooth] plot (\x,{sqrt(max(0,2*\x*(\x+2)*(\x+9)))});
		\node[pt,label={above left:$P$}] at (-3,6) {};
		\node[pt,label={left:$2P$}] at (1.125,8.4375) {};
		\node[pt,label={[label distance=1pt]right:$3P$}] at (-8.9504,-2.4838) {};
		\node[pt,label={below:$5P$}] at (-2.5871,-4.4137) {};
		\node[pt,label={right:$-5P$}] at (-2.5871,4.4137) {};
		\node[pt,label={above left:$T$}] at (-9,0) {};
		\node[below] at (-2.5,-14) {\small window};
		\node at (-5.8,9.5) {$\Omega$};
		\node at (2.3,11) {$\Omega_\infty$};
		\node[below,gray] at (-1.75,-0.2) {\small $-2$};
		\node[below left,gray] at (0,-0.2) {\small $0$};
	\end{tikzpicture}
	\caption{The real points of $E_{2,1}:Y^2=2w(w+2)(w+9)$: the oval $\Omega$ and the unbounded branch $\Omega_\infty$. The window is $-3<w<-2$, and the thick arc over it carries the rational $\Delta$-triples for perfect squares with $z/y=2$. The point $P=(-3,6)$ sits at the left end of the window (\Cref{rem:ell-ends}). Its multiples alternate between the two pieces (\Cref{cor:ell-parity}); $5P$ and its mirror image $-5P$ lie over the window. Here $T=(-9,0)$, and the oval meets the $w$-axis at $w=-9$ and $w=-2$.}
	\label{fig:ell-curve}
\end{figure}

Which multiples of a point land in the window? To answer this, we look at the real picture. Since $n\,ab>0$, the sign of $f(w)=n\,ab\,w(w+ab)(w+s^2)$ is
\[
\begin{array}{c|c|c|c|c}
	w & w<-s^2 & -s^2\leq w\leq-ab & -ab<w<0 & w\geq0\\\hline
	f(w) & <0 & \geq0 & <0 & \geq0
\end{array}
\]
As $Y^2=f(w)\geq0$, the real points of $E$ form two pieces (\Cref{fig:ell-curve}):
\begin{itemize}
	\item the \emph{oval} $\Omega$, made of the real points with $-s^2\leq w\leq-ab$; it is a closed curve;
	\item the \emph{unbounded branch} $\Omega_\infty$, made of the real points with $w\geq0$, together with $O$.
\end{itemize}
A real point of $E$ with $w<0$ lies on $\Omega$. In particular, the rational point $T=(-s^2,0)$ lies on $\Omega$, and $2T=O$ because $T=-T$. Since $bs<s^2$, the window lies inside the $w$-range of the oval.

We write $\mathbb{R}/\mathbb{Z}$ for the quotient of the additive group $\mathbb{R}$ by its subgroup $\mathbb{Z}$, and $\overline{x}$ for the class of $x\in\mathbb{R}$. Thus $\overline{x}=\overline{y}$ if and only if $x-y\in\mathbb{Z}$, and $\overline{x}+\overline{y}=\overline{x+y}$. A subset $\mathcal{U}\subseteq\mathbb{R}/\mathbb{Z}$ is \emph{open} if $\{x\in\mathbb{R}:\overline{x}\in \mathcal{U}\}$ is open in $\mathbb{R}$. We give $E(\mathbb{R})$ its usual topology as a curve in the real projective plane; away from $O$, it is the topology induced by $\mathbb{R}^2$.

\begin{theorem}[Real points; see {\cite[Chapter~II]{silverman-tate}} and {\cite[Chapter~V, \S2]{silverman-advanced}}]
	\label{fact:ell-real}
	The addition $E(\mathbb{R})\times E(\mathbb{R})\to E(\mathbb{R})$ is continuous, and $\Omega_\infty$ and $\Omega$ are the two connected components of $E(\mathbb{R})$. The component $\Omega_\infty$ is a subgroup of $E(\mathbb{R})$ of index $2$, so $\Omega$ is its non-trivial coset. Moreover, there is a bijection $\phi\colon\Omega_\infty\to\mathbb{R}/\mathbb{Z}$ which is a homeomorphism and satisfies $\phi(Q+Q')=\phi(Q)+\phi(Q')$ for all $Q,Q'\in\Omega_\infty$.
\end{theorem}

In other words, $\Omega_\infty$ is a copy of $\mathbb{R}/\mathbb{Z}$, both as a group and as a topological space. Note that $\phi(O)=\overline{0}$, because $\phi(O)=\phi(O+O)=\phi(O)+\phi(O)$; and $\phi(kQ)=k\,\phi(Q)$ for $Q\in\Omega_\infty$ and $k\geq1$.

\begin{corollary}[Parity]
	\label{cor:ell-parity}
	If $P_1,P_2\in\Omega$, then $P_1+P_2\in\Omega_\infty$; if $P_1\in\Omega$ and $P_2\in\Omega_\infty$, then $P_1+P_2\in\Omega$. In particular, if $P\in\Omega$ and $j\geq1$, then $jP$ lies on $\Omega$ if $j$ is odd and on $\Omega_\infty$ if $j$ is even.
\end{corollary}

\begin{proof}
	The quotient $E(\mathbb{R})/\Omega_\infty$ has order $2$ by \Cref{fact:ell-real}.
\end{proof}

So the odd multiples of a point of the oval stay on the oval, and they are the only multiples that can produce triples. To see that they actually enter the window, we use Kronecker's theorem.

\begin{theorem}[Kronecker; see {\cite[Chapter~XXIII, Theorem~438]{hardy-wright}}]
	\label{fact:kronecker}
	Let $\vartheta$ be an irrational number, let $\alpha\in\mathbb{R}$, $\varepsilon>0$ and $N\in\mathbb{N}$. Then there are integers $k>N$ and $p$ such that $|k\vartheta-p-\alpha|<\varepsilon$.
\end{theorem}

\begin{corollary}
	\label{cor:kronecker}
	Let $\vartheta$ be irrational and let $\mathcal{U}\subseteq\mathbb{R}/\mathbb{Z}$ be a non-empty open set. Then $\overline{k\vartheta}\in \mathcal{U}$ for infinitely many integers $k\geq1$.
\end{corollary}

\begin{proof}
	Pick $\overline{\alpha}\in \mathcal{U}$. The set $\{x\in\mathbb{R}:\overline{x}\in \mathcal{U}\}$ is open and contains $\alpha$, so it contains an interval $(\alpha-\varepsilon,\alpha+\varepsilon)$ with $\varepsilon>0$. For every $N\in\mathbb{N}$, \Cref{fact:kronecker} gives integers $k>N$ and $p$ with $k\vartheta-p\in(\alpha-\varepsilon,\alpha+\varepsilon)$, hence $\overline{k\vartheta}=\overline{k\vartheta-p}\in \mathcal{U}$.
\end{proof}

\begin{lemma}[Density]
	\label{lem:ell-density}
	Let $P\in\Omega\cap E(\mathbb{Q})$ have infinite order, and let $J$ be a non-empty open interval contained in $(-s^2,-ab)$. Then $w(jP)\in J$ for infinitely many odd $j\geq1$.
\end{lemma}

\begin{proof}
	Let $\tau\colon E(\mathbb{R})\to E(\mathbb{R})$ be the translation $\tau(Q)=P+Q$. It is continuous by \Cref{fact:ell-real}, and bijective, with inverse $Q\mapsto(-P)+Q$.
	
	\emph{Step 1: a target set on the branch.} Let $U$ be the set of points $Q\neq O$ of $E(\mathbb{R})$ with $w(Q)\in J$. It is open in $E(\mathbb{R})$, because $w$ is continuous away from $O$, and it is non-empty: for $w_0\in J$ we have $f(w_0)>0$, so $\bigl(w_0,\sqrt{f(w_0)}\bigr)\in U$. Since its points have $w<0$, $U\subseteq\Omega$. Put $V=\tau^{-1}(U)$, a non-empty open subset of $E(\mathbb{R})$. If $Q\in V\cap\Omega$, then $P+Q\in\Omega_\infty$ by \Cref{cor:ell-parity}, contradicting $P+Q\in U\subseteq\Omega$. Hence $V\subseteq\Omega_\infty$, and $V$ is open in $\Omega_\infty$.
	
	\emph{Step 2: even multiples enter $V$.} By \Cref{cor:ell-parity}, $2P\in\Omega_\infty$. Choose $\vartheta\in\mathbb{R}$ with $\phi(2P)=\overline{\vartheta}$; then $\phi(2kP)=k\,\overline{\vartheta}=\overline{k\vartheta}$ for $k\geq1$. The number $\vartheta$ is irrational: if $\vartheta=p/q$ with $q\geq1$, then $\phi(2qP)=\overline{p}=\overline{0}=\phi(O)$, so $2qP=O$ because $\phi$ is injective, contradicting the fact that $P$ has infinite order. Since $\phi$ is a homeomorphism, $\phi(V)$ is a non-empty open subset of $\mathbb{R}/\mathbb{Z}$, and \Cref{cor:kronecker} gives infinitely many $k\geq1$ with $\overline{k\vartheta}\in\phi(V)$, that is, $2kP\in V$.
	
	\emph{Step 3: conclusion.} For each such $k$, the point $(2k+1)P=\tau(2kP)$ lies in $U$, that is, $w\bigl((2k+1)P\bigr)\in J$.
\end{proof}

\begin{ejemplo}
	\label{ex:ell-24-window}
	We continue \Cref{ex:ell-24-chords}. The odd multiples of $P_0$ on $E^{(24)}_{2,1}$ have the following $w$-coordinates (window: $-3<w<-2$):
	\[
	\begin{array}{c|cccccc}
		j & 1 & 3 & 5 & 7 & 9 & 11\\\hline
		w(jP_0)\approx & -2.250 & -4.622 & -8.632 & -7.169 & -3.295 & -2.008\\
		\text{in the window?} & \text{yes} & \text{no} & \text{no} & \text{no} & \text{no} & \text{yes}
	\end{array}
	\]
	The point $11P_0$ (or its mirror image) gives, through \eqref{eq:ell-backward}, a new rational $\Delta$-triple for $24$,
	\[
	(x,y,z)\approx(0.16324,\ 0.24391,\ 0.48783),\qquad z/y=2,
	\]
	which satisfies \eqref{eq:tool-delta} exactly. The number $L$ of \Cref{lem:ell-rational-integer}(3) has $173$ decimal digits, and $(Lx,Ly,Lz)$ is a $\Delta$-triple for $24L^2$ that is not proportional to $(3,4,8)$. Among the first twenty odd multiples of $P_0$, seven lie in the window.
\end{ejemplo}

\section{Squares: duplicated versus generic triples}
\label{subsec:ell-squares}

For perfect squares we may take $n=1$ (\Cref{rem:ell-scaling}), and the curve is
\[
E_{a,b}:\qquad Y^2=ab\,w\,(w+ab)\,(w+s^2),
\]
with Weierstrass model
\[
V^2=u\,(u+a^2b^2)\,(u+ab\,s^2),
\]
where $u=ab\,w$ and $V=ab\,Y$. By \Cref{rem:ell-ends}, $E_{a,b}$ always contains the rational point
\[
P=\bigl(-bs,\ ab^2s\bigr),
\]
coming from the degenerate triple $(b/a,\,b/a,\,1)$. The point $P$ lies on the oval, at the left end of the window. It is not itself a $\Delta$-triple, but by \Cref{lem:ell-density}, if it has infinite order its odd multiples fill the oval, and in particular enter the window. Everything thus depends on the order of $P$.

\begin{lemma}
	\label{lem:ell-order-P}
	Let $a\geq b\geq1$ be coprime, $s=a+b$, and $P=(-bs,ab^2s)\in E_{a,b}(\mathbb{Q})$.
	\begin{enumerate}
		\item In the Weierstrass model, $u(2P)=\dfrac{(a^2-b^2)^2}{4}$.
		\item If $a=b=1$, then $P$ has order $4$.
		\item If $a>b$, then $P$ has infinite order.
	\end{enumerate}
\end{lemma}

\begin{proof}
	In the Weierstrass model $V^2=u^3+Au^2+Bu$ we have $B=a^2b^2\cdot ab\,s^2=a^3b^3s^2$, and $P$ becomes $(u_0,V_0)=(-ab^2s,\ a^2b^3s)$.
	
	(1) We have $u_0^2-B=a^2b^4s^2-a^3b^3s^2=a^2b^3s^2(b-a)$ and $4V_0^2=4a^4b^6s^2$. By \Cref{lem:ell-duplication},
	\[
	u(2P)=\frac{a^4b^6s^4(a-b)^2}{4a^4b^6s^2}=\frac{s^2(a-b)^2}{4}=\frac{(a^2-b^2)^2}{4}.
	\]
	
	(2) If $a=b=1$, then $u(2P)=0$, so $2P=(0,0)$. The point $(0,0)$ is its own mirror image, so $(0,0)+(0,0)=(0,0)+\bigl(-(0,0)\bigr)=O$, that is, $2(2P)=O$. Hence $4P=O$, while $P\neq O$ and $2P\neq O$: the order of $P$ is $4$. Explicitly, the multiples of $P=(-2,2)$ are $P$, $2P=(0,0)$, $3P=(-2,-2)$ and $4P=O$.
	
	(3) Let $a>b$. By \Cref{fact:ell-NL}, it suffices to show that $u(2P)$ or $u(4P)$ is not an integer. Since $a,b$ are coprime, they are not both even.
	
	\emph{Case 1: $a$ and $b$ have different parity.} Then $a^2-b^2$ is odd, and by (1) $u(2P)$ is an odd number divided by $4$, which is not an integer.
	
	\emph{Case 2: $a$ and $b$ are both odd.} Then
	\[
	a+b=2g,\qquad a-b=2h,\qquad a^2+b^2=2e
	\]
	with positive integers $g,h,e$, and $e$ is odd because $a^2+b^2\equiv2\pmod4$. Now $u_2:=u(2P)=4g^2h^2$ is an integer, so we double once more. Since $(a^2-b^2)^2+4a^2b^2=(a^2+b^2)^2$ and $(a-b)^2+4ab=(a+b)^2$, we get
	\begin{align*}
		u_2+a^2b^2&=\frac{(a^2-b^2)^2+4a^2b^2}{4}=\frac{(a^2+b^2)^2}{4}=e^2,\\
		u_2+ab\,s^2&=\frac{(a+b)^2(a-b)^2+4ab(a+b)^2}{4}=\frac{(a+b)^4}{4}=4g^4 .
	\end{align*}
	Hence $V(2P)^2=u_2\,(u_2+a^2b^2)(u_2+ab\,s^2)=16\,g^6h^2e^2\neq0$. Moreover $u_2^2-B=16g^4h^4-4g^2a^3b^3=4g^2\,(4g^2h^4-a^3b^3)$. By \Cref{lem:ell-duplication},
	\[
	u(4P)=\frac{16g^4\,(4g^2h^4-a^3b^3)^2}{4\cdot16\,g^6h^2e^2}=\frac{(4g^2h^4-a^3b^3)^2}{4\,g^2h^2e^2}.
	\]
	The numerator is odd, since $a^3b^3$ is odd and $4g^2h^4$ is even, while the denominator is even. So $u(4P)$ is not an integer.
	
	In both cases $P$ has infinite order.
\end{proof}

The formula in (1) contains the whole dichotomy: $u(2P)=0$ exactly when $a=b$. In the duplicated case $2P$ is the point $(0,0)$ of order $2$; in the generic case $2P$ is a point of infinite order.

\begin{theorem}
	\label{thm:ell-every-quotient}
	Let $a>b\geq1$ be coprime integers. Then there are infinitely many pairwise non-proportional $\Delta$-triples $(x,y,z)$ of perfect squares with $\frac zy=\frac ab$.
\end{theorem}

\begin{proof}
	Let $s=a+b$ and $P=(-bs,ab^2s)\in E_{a,b}(\mathbb{Q})$. Since $-s^2\leq-bs\leq-ab$, the point $P$ lies on the oval, and by \Cref{lem:ell-order-P}(3) it has infinite order. By \Cref{lem:ell-density}, applied to the window $J=(-bs,-ab)$, infinitely many odd multiples $jP$ have $w(jP)$ in the window. These multiples are pairwise distinct, and at most two points of $E$ share the same $w$, so we obtain infinitely many distinct values $w$ in the window. For each of them, the point $(w,|Y|)$ gives an element of $\mathcal{R}_1(a/b)$ by \Cref{prop:ell-curve}, and distinct values of $w$ give distinct elements, because $x/y=-w/(bs)$ by \eqref{eq:ell-forward}. Hence $\mathcal{R}_1(a/b)$ is infinite, and \Cref{cor:ell-classes} gives infinitely many pairwise non-proportional $\Delta$-triples of perfect squares with $z/y=a/b$.
\end{proof}

\begin{remark}[The duplicated case on the curve]
	\label{rem:ell-duplicated}
	For $a=b=1$ the curve is $E_{1,1}:Y^2=w(w+1)(w+4)$, the window is $-2<w<-1$, and $P=(-2,2)$ has order $4$ (\Cref{lem:ell-order-P}(2)). Its odd multiples are $P$ and $3P=-P$, both with $w=-2$: they stay at the end of the window, which corresponds to the degenerate triple $(1,1,1)$. So the mechanism of \Cref{thm:ell-every-quotient} produces nothing.
	
	By \Cref{prop:ell-curve} and \Cref{lem:ell-rational-integer}(2),(3), a perfect square has a duplicated $\Delta$-triple if and only if $E_{1,1}$ has a rational point over the window. The curve $E_{1,1}$ has conductor $24$ and rank $0$ (a $2$-descent, e.g.\ with \texttt{ellrank} of PARI/GP \cite{pari}, shows that the rank is $0$), and its rational points are the eight points of finite order
	\[
	O,\quad(0,0),\quad(-1,0),\quad(-4,0),\quad(-2,\pm2),\quad(2,\pm6),
	\]
	none of which has $-2<w<-1$. This gives an alternative proof of \cite[Theorem~2.6]{gs-squares}, based on the rank of $E_{1,1}$ instead of the theorem of Euler on products of four terms in arithmetic progression (see \cite{hlst}).
\end{remark}

\begin{remark}[Other points of the oval]
	\label{rem:ell-other-points}
	By \Cref{cor:ell-parity}, the points $2iP+T$ also lie on the oval, and they also produce triples when they fall in the window. For example, the least perfect square with a $\Delta$-triple with $z/y=2$ is $180^2$, with the triple $(75,81,162)$. Its point on $E_{2,1}$ is $\bigl(-\tfrac{25}{9},\tfrac{140}{27}\bigr)=-2P+T$, not an odd multiple of $P$.
\end{remark}

\begin{ejemplo}
	\label{ex:ell-squares}
	\begin{enumerate}
		\item For $z/y=2$ we have $E_{2,1}:Y^2=2w(w+2)(w+9)$, $P=(-3,6)$ and $2P=(\tfrac98,\tfrac{135}{16})$; in the Weierstrass model $u(2P)=\tfrac94\notin\mathbb{Z}$, as in Case~1 of the proof of \Cref{lem:ell-order-P}. The first odd multiple of $P$ over the window is $5P$ (see \Cref{fig:ell-curve}). Through \eqref{eq:ell-backward} and \Cref{lem:ell-rational-integer}(3) it gives the $\Delta$-triple
		\[
		(56555348268241,\ 65581507294081,\ 131163014588162)
		\]
		for $N^2$, where $N=164343419150038$.
		\item For $z/y=3$ and $z/y=3/2$, the least perfect squares with such a $\Delta$-triple are $18837^2$ and $3900^2$ (\Cref{sec:conclusion}). On the curve, the first odd multiples of $P$ over the window are $7P$ and $5P$, respectively, and they give $\Delta$-triples for $N^2$ with
		\[
		N=488650465808508899094963
		\]
		and
		\[
		N=3584529757110101799486 .
		\]
		\item Small examples with other quotients are $30^2$ with $(9,10,25)$, where $z/y=5/2$, and $60^2$ with $(9,16,18)$, where $z/y=9/8$.
	\end{enumerate}
\end{ejemplo}

\section{Arbitrary $n$: the rank decides}
\label{subsec:ell-rank}

For $n$ not a perfect square, the degenerate point of \Cref{rem:ell-ends} is not rational, and there is no free point to start from. The question becomes whether $E^{(n)}_{a,b}$ has infinitely many rational points. Recall the following theorem.

\begin{theorem}[Mordell; see {\cite[Chapter~III]{silverman-tate}}]
	\label{fact:ell-mordell}
	The group $E(\mathbb{Q})$ is finitely generated. Hence $E(\mathbb{Q})\cong\mathbb{Z}^\rho\oplus F$ with $F$ finite; the integer $\rho\geq0$ is the \emph{rank} of $E$.
\end{theorem}

In particular, $E(\mathbb{Q})$ is infinite if and only if $\rho\geq1$, if and only if $E(\mathbb{Q})$ has a point of infinite order.

\begin{theorem}
	\label{prop:ell-rank}
	Let $n\in\mathbb{N}$ and let $a\geq b\geq1$ be coprime. Then $\mathcal{R}_n(a/b)$ is infinite if and only if $E^{(n)}_{a,b}$ has positive rank. In that case there are infinitely many pairwise non-proportional $\Delta$-triples with $z/y=a/b$ among the numbers $\kappa^2n$, $\kappa\in\mathbb{N}$.
\end{theorem}

\begin{proof}
	If the rank is $0$, then $E^{(n)}_{a,b}(\mathbb{Q})$ is finite, and so is $\mathcal{R}_n(a/b)$ by \Cref{prop:ell-curve}.
	
	Assume the rank is positive, and let $Q$ be a rational point of infinite order. If $Q\in\Omega_\infty$, replace $Q$ by $Q+T$. Then $Q+T\in\Omega$ by \Cref{cor:ell-parity}, and $Q+T$ still has infinite order, because $2(Q+T)=2Q$. So we may assume $Q\in\Omega$. As in the proof of \Cref{thm:ell-every-quotient}, \Cref{lem:ell-density} and \Cref{prop:ell-curve} show that $\mathcal{R}_n(a/b)$ is infinite. The last assertion follows from \Cref{cor:ell-classes}.
\end{proof}

\begin{corollary}
	\label{cor:ell-triple-infinite-order}
	Let $(x,y,z)$ be a $\Delta$-triple for $n$, with $z/y=a/b$ in lowest terms. If the point \eqref{eq:ell-forward} of $E^{(n)}_{a,b}$ has infinite order, then the numbers $\kappa^2n$, $\kappa\in\mathbb{N}$, have infinitely many pairwise non-proportional $\Delta$-triples with $z/y=a/b$.
\end{corollary}

\begin{proof}
	The point lies over the window, hence on the oval, and has infinite order. Proceed as in the proof of \Cref{thm:ell-every-quotient}.
\end{proof}

For instance, $n=24$ and the triple $(1,\tfrac43,\tfrac83)$, that is, the $\Delta$-triple $(3,4,8)$ for $216$, satisfy the hypothesis by \Cref{ex:ell-24-chords}. The hypothesis was tested with PARI/GP \cite{pari} on every $\Delta$-triple of every $n\leq10^5$, as follows. Each such triple $(x,y,z)$ has its own quotient $z/y=a/b$ in lowest terms, hence its own curve $E^{(n)}_{a,b}$, on which \eqref{eq:ell-forward} places it as a rational point; the order of that point was computed in the group $E^{(n)}_{a,b}(\mathbb{Q})$, using the Weierstrass model \eqref{eq:ell-weierstrass}. There are $4058$ such triples, for $1780$ values of $n$, and in every case the point turned out to have infinite order. So, by \Cref{cor:ell-triple-infinite-order}, each of these $1780$ numbers satisfies the conclusion of \Cref{prop:ell-rank} for the corresponding quotients $a/b$.

\begin{problem}
	\label{prob:ell-rank}
	For fixed coprime $a\geq b\geq1$, determine the squarefree integers $n$ for which $E^{(n)}_{a,b}$ has positive rank.
\end{problem}

Since $E^{(n)}_{a,b}$ is the curve $n(Y')^2=ab\,w(w+ab)(w+s^2)$ (put $Y'=Y/n$), that is, the quadratic twist of $E_{a,b}$ by $n$, \Cref{prob:ell-rank} asks which quadratic twists of the fixed curve $E_{a,b}$ have positive rank. In fact, \Cref{prob:ell-rank} is a $\theta$-congruent number problem in the sense of Fujiwara \cite{fujiwara}. Write $\cos\theta=(a^2+b^2)/s^2=\gamma/\delta$ in lowest terms. The substitution $x=nab\,(w+ab)$, $y=nab\,Y$ turns \eqref{eq:ell-curve} into $y^2=x\,(x-na^2b^2)\bigl(x+nab\,(a^2+ab+b^2)\bigr)$, and since $(\delta-\gamma)/(\delta+\gamma)=ab/(a^2+ab+b^2)$, the change $x\mapsto(\delta-\gamma)^2x/(ab)^2$ makes it the curve
\[
\mathcal{E}_{m,\theta}:\quad y^2=x\bigl(x+(\delta+\gamma)m\bigr)\bigl(x-(\delta-\gamma)m\bigr),\qquad m=n(\delta-\gamma).
\]
A triangle with rational sides $A,B,C$, an angle $\theta$ between $A$ and $B$, and area $m\sqrt{\delta^2-\gamma^2}$, that is, with $AB=2\delta m$, gives the point $\bigl(C^2/4,\,C(A^2-B^2)/8\bigr)$ of $\mathcal{E}_{m,\theta}$, by the law of cosines. Conversely, if $Q$ is a rational point with $2Q\neq O$, then $x(2Q)$, $x(2Q)+(\delta+\gamma)m$ and $x(2Q)-(\delta-\gamma)m$ are the squares $C^2/4$, $(A+B)^2/4$ and $(A-B)^2/4$ for such a triangle. By \cite{fujiwara}, provided the squarefree part of $m$ is not $1$, $2$, $3$ or $6$, such a triangle exists exactly when $\mathcal{E}_{m,\theta}$ has positive rank.

Since $\cos\theta\geq1/2$, the angle $\theta$ lies in $(0,\pi/3]$, and the classical congruent number problem ($\theta=\pi/2$) never occurs. For $a=b$ we get $\theta=\pi/3$, and $E^{(n)}_{1,1}$ becomes the $\pi/3$-congruent number curve $y^2=x(x-n)(x+3n)$. For $n=1$, the point $P=(-2,2)$ of \Cref{lem:ell-order-P}(2) becomes $(-1,2)$, and $2P=(1,0)$ is the point given by the equilateral triangle of side $2$: the degenerate point $P$ is a half of the equilateral triangle. For $a/b=5/2$ we have $\cos\theta=29/49$ and $m=20$; doubling the point of $\mathcal{E}_{20,\theta}$ given by the $\Delta$-triple $(9,10,25)$ for $30^2$ yields the triangle with sides $84$, $70/3$, $218/3$, an angle $\theta$ between the first two, and area $20\sqrt{49^2-29^2}$.

\section{What the method does not see}
\label{subsec:ell-limits}

By \Cref{lem:ell-rational-integer}, a rational $\Delta$-triple for $n$ describes $\Delta$-triples for all the numbers $\kappa^2n$ with $\kappa$ a multiple of a certain $L$. Thus the curves $E^{(n)}_{a,b}$ see a number and all its multiples by squares at the same time, and they cannot tell which of these numbers is the smallest one carrying a given triple. In particular, the method gives no information on $\Delta$-primitivity. For instance, \Cref{thm:ell-every-quotient} does not claim that any of the squares it produces is $\Delta$-primitive, and the square $180^2$ of \Cref{rem:ell-other-points} is not, since $30\in\mathcal{S}$ divides $180$ (\Cref{sec:prelim}). Questions on $\Delta$-primitive numbers require arithmetic arguments of a different kind, such as those of \cite{gs-squares}.

\section{Concluding remarks and future work}
\label{sec:conclusion}

Fixing the quotient $z/y=a/b$ turns the search for $\Delta$-triples into the study of rational points on the explicit elliptic curves $E^{(n)}_{a,b}$, on which the chord--tangent construction produces new $\Delta$-triples from old ones. For perfect squares the picture is complete: everything is decided by the order of the rational point $P$ coming from the degenerate triple, which is finite exactly in the duplicated case $a=b$. For arbitrary $n$ the picture is decided by the rank of $E^{(n)}_{a,b}$, which leads to \Cref{prob:ell-rank}. We end with two further problems.

\emph{The size of the least square.} For coprime $a>b\geq1$, let $\nu(a/b)$ be the least $N\in\mathbb{N}$ such that $N^2$ has a $\Delta$-triple with $z/y=a/b$; it exists by \Cref{thm:ell-every-quotient}. We have $\nu(2)=180$ (\Cref{rem:ell-other-points}), and a direct search gives $\nu(3/2)=3900$, with the $\Delta$-triple $(845,1250,1875)$, and $\nu(3)=18837$, with the $\Delta$-triple $(3969,4563,13689)$. These are much smaller than the squares produced in \Cref{ex:ell-squares} by the first odd multiples of $P$ over the window.

\begin{problem}
	\label{prob:least-square}
	Estimate $\nu(a/b)$ in terms of $a$ and $b$.
\end{problem}

The natural tool for \Cref{prob:least-square} is the theory of heights on elliptic curves, which measures how fast the size of the coordinates of $jP$ grows with $j$.

\emph{Quotients of $\Delta$-primitive squares.} By \Cref{subsec:ell-limits}, the curves cannot decide whether the squares they produce are $\Delta$-primitive. Some quotients do occur on $\Delta$-primitive squares: $30^2$ is $\Delta$-primitive (\Cref{sec:prelim}) and has the $\Delta$-triple $(9,10,25)$, with $z/y=5/2$. Others seem not to: a computer search shows that, for $N\leq10^4$, the integers $N$ such that $N^2$ has a $\Delta$-triple with $z/y=2$ are exactly the multiples of $180$ in that range, and none of these squares is $\Delta$-primitive.

\begin{problem}
	\label{prob:primitive-quotients}
	Determine the coprime pairs $a>b\geq1$ for which some $\Delta$-primitive square has a $\Delta$-triple with $z/y=a/b$. In particular, is there a $\Delta$-primitive square with a $\Delta$-triple with $z/y=2$?
\end{problem}

\section*{Acknowledgements}
This work was partially supported by: Universidad Nacional de San Luis (UNSL), grants PI UNSL 03-1326, PROIPRO 03-2923, and PROINI 03-124; Consejo Nacional de Investigaciones Cient\'ificas y T\'ecnicas (CONICET), grant PIP 11220220100068CO.


\begin{thebibliography}{99}
	\bibitem{erdosproblems}
	T.F.~Bloom, \emph{Erd\H{o}s problems}, \url{https://www.erdosproblems.com}, Problem~885, accessed November 2026.

	\bibitem{bremner}
	A.~Bremner, \emph{On a problem of Erd\H{o}s related to common factor differences}, Int.\ J.\ Number Theory \textbf{15} (2019), no.~5, \href{https://doi.org/10.1142/S1793042119500581}{doi:10.1142/S1793042119500581}.

	\bibitem{er-factor-difference}
	P.~Erd\H{o}s and M.~Rosenfeld, \emph{The factor-difference set of integers}, Acta Arith.\ \textbf{79} (1997), no.~4, 353--359.

    \bibitem{fujiwara}
    M.~Fujiwara, \emph{$\theta$-congruent numbers}, in: Number Theory (Eger, 1996), de Gruyter, Berlin, 1998, 235--241.

	\bibitem{gs-squares}
	F.~Golfieri and V.N.~Schv\"ollner, \emph{On the $\Delta$-primitive squares conjecture}, arXiv:2611.11574, 2026.

	\bibitem{hardy-wright}
	G.H.~Hardy and E.M.~Wright, \emph{An Introduction to the Theory of Numbers}, 6th ed., revised by D.R.~Heath-Brown and J.H.~Silverman, Oxford University Press, Oxford, 2008.

	\bibitem{hlst}
	N.~Hirata-Kohno, S.~Laishram, T.N.~Shorey and R.~Tijdeman, \emph{An extension of a theorem of Euler}, Acta Arith.\ \textbf{129} (2007), 71--102.

	\bibitem{jimenez-urroz}
	J.~Jim\'enez-Urroz, \emph{A note on a conjecture of Erd\H{o}s and Rosenfeld}, J.\ Number Theory \textbf{78} (1999), no.~1, 140--143.

	\bibitem{pari}
	The PARI~Group, \emph{PARI/GP}, version 2.17.2, Univ.\ Bordeaux, \url{https://pari.math.u-bordeaux.fr/}.

	\bibitem{vnsigma.factor.graph}
	A.~Pastine and V.N.~Schv\"ollner, \emph{Simple factor graphs associated with split graphs}, arXiv:2512.24252, 2025.

	\bibitem{pomerance}
	C.~Pomerance, \emph{A tale of two sieves}, Notices Amer.\ Math.\ Soc.\ \textbf{43} (1996), no.~12, 1473--1485.

	\bibitem{vnsigma.delta.property}
	V.N.~Schv\"ollner, \emph{The $\Delta$ property: a bridge between split graphs and Number Theory}, arXiv:2605.25264, 2026.

	\bibitem{vnsigma.factor.graph.paths}
	V.N.~Schv\"ollner and A.~Pastine, \emph{Induced paths and cycles in factor graphs of split graphs}, arXiv:2603.14061, 2026.

	\bibitem{silverman-advanced}
	J.H.~Silverman, \emph{Advanced Topics in the Arithmetic of Elliptic Curves}, Graduate Texts in Mathematics \textbf{151}, Springer, New York, 1994.

	\bibitem{silverman-tate}
	J.H.~Silverman and J.T.~Tate, \emph{Rational Points on Elliptic Curves}, 2nd ed., Undergraduate Texts in Mathematics, Springer, Cham, 2015.

	\bibitem{tao-vu}
	T.~Tao and V.~Vu, \emph{Sum-free sets in groups: a survey}, arXiv:1603.03071, 2016.

\end{thebibliography}
\end{document}